\documentclass{article}[12pt]
\usepackage{bm, amsmath, amsfonts, amssymb, cite}
\usepackage{graphicx}

\makeatletter \@addtoreset{equation}{section} \makeatother

\newtheorem{theo}{Theorem}[section]

\newtheorem{coro}[theo]{Corollary}
\newtheorem{defi}[theo]{Definition}
\newtheorem{lem}[theo]{Lemma}

\def\F{\mathbb F}
\def\Z{\mathbb Z}
\def\N{\mathbb N}

\def\a{\alpha}

\def\e{{\boldsymbol \varepsilon}}
\def\den{{\rm denom}}
\def\num{{\rm numer}}
\def\qdis{{\rm qdis}}

\def\mid{{\,|\,}}

\def\pf{\noindent {\it Proof.\ }}
\def\qed{\hfill \rule{4pt}{7pt}}
\begin{document}

\begin{center}
{\Large Applicability of the $q$-Analogue of Zeilberger's
Algorithm}

\vskip 3mm

 William Y.C. Chen$^1$,
Qing-Hu Hou$^2$ and Yan-Ping Mu$^3$\\[5pt]
Center for Combinatorics, LPMC \\
Nankai University, Tianjin 300071, P. R. China

\vskip 3mm

 E-mail: $^1$chen@nankai.edu.cn,
$^2$hou@nankai.edu.cn, $^3$myphb@eyou.com

\end{center}

\begin{abstract}

The applicability or terminating condition for the ordinary case
of Zeilberger's algorithm was recently obtained by Abramov. For
the $q$-analogue, the question of whether a bivariate
$q$-hypergeometric term has a $qZ$-pair remains open. Le has found
a solution to this problem when the given bivariate
$q$-hypergeometric term is a rational function in certain powers
of $q$. We solve the problem for the general case by giving a
characterization of bivariate $q$-hypergeometric terms for which
the $q$-analogue of Zeilberger's algorithm terminates. Moreover,
we give an algorithm to determine whether a bivariate
$q$-hypergeometric term has a $qZ$-pair.
\end{abstract}



\section{Introduction}
Zeilberger's algorithm \cite{Gra,Pet-W-Z,Zeil}, also known as the method of
{\it creative telescoping}, is devised for proving hypergeometric identities of
the form
$$
\sum_{k=-\infty}^{\infty} F(n,k) = f(n),
$$
where $F(n,k)$ is a bivariate hypergeometric term and $f(n)$ is a
given function (for most cases a hypergeometric term plus a
constant). The algorithm can be easily adapted to the $q$-case,
which is called the $q$-analogue of Zeilberger's algorithm
\cite{Koo,Bo-Ko,P-R,W-Z}. Let $N$ and $K$ be the shift operators
with respect to $n$ and $k$ respectively, defined by
\[ N T(n,k)=T(n+1,k) \qquad \mbox{and} \qquad K T(n,k)=T(n,k+1).\]
Given a bivariate $q$-hypergeometric term $T(n,k)$, the
$q$-analogue of Zeilberger's algorithm  aims to find a {\it
$qZ$-pair} $(L,G)$, where  $L$ is a linear difference operator
with coefficients in the ring of polynomials in $q^n$
$$
L=a_0(q^n) N^0 + a_1(q^n) N^1 +\cdots+a_r(q^n) N^r
$$
and $G$ is a bivariate $q$-hypergeometric term $G(n,k)$ such that
$$
L T(n,k)=(K-1)G(n,k).
$$
Zeilberger's algorithm has been widely used as a powerful tool to
prove hypergeometric identities. It was an open question when the
algorithm terminates. This problem was solved recently by Abramov
\cite{Abra-1,Abra-2}. For the $q$-analogue of Zeilberger's
algorithm, Le \cite{Le} found a solution to the termination
problem  for the case of rational functions. In this paper we
provide a complete solution for the general $q$-case.

We begin with an additive decomposition of univariate
$q$-hypergeometric terms. Using this decomposition, a univariate
$q$-hypergeometric term $T(n)$ can be represented as
\[ T(n)=(N-1)T_1(n) + T_2(n),\]
 where $T_1(n)$ and $T_2(n)$ are $q$-hypergeometric terms, and
 $T_2(n)$ has the following form
$$T_2(n)={u_1(q^n) \over u_2(q^n)} \prod_{j=n_0}^{n-1} {f_1(q^j) \over f_2(q^j)},$$
where $u_1,u_2,f_1,f_2$ are polynomials and for any integer $m$, $u_2(x)$ and
$u_2(xq^m)$ have no common factors except for a power of $x$. Consequently, a
bivariate $q$-hypergeometric term $T(n,k)$ can be decomposed as
\begin{equation}
\label{ddd} T(n,k)=(K-1)T_1(n,k) + T_2(n,k)
\end{equation}
such that
$$T_2(n,k)=T(n,k_0) V(q^n,q^k) \prod_{j=k_0}^{k-1}F(q^n,q^j),$$
where $V,F$ are rational functions and the denominator $v_2$ of $V$ satisfies
the conditions that for any integer $m$, $v_2(x,y)$ and $v_2(x,yq^m)$ have no
common factors except for a power of $y$. The polynomial $v_2(x,y)$ with the
above property is called $\e_y$-free. We should note that the above
decomposition does not solve the minimal additive decomposition problem and is
not unique (see \cite{Ab-Pet} for a precise definition). However, for the
purpose of constructing a $qZ$-pair, it turns out that one may choose any
decomposition.

Then we consider the structure of bivariate $q$-hypergeometric
terms. The structure of ordinary hypergeometric terms has been
studied by Ore \cite{Ore}, Sato-Shintani-Muro\cite{Sato},
Abramov-Petkov\v{s}ek \cite{Ab-Pet-1} and Hou \cite{Hou}. To a
large extent, the $q$-case is analogous to the ordinary case. For
each bivariate $q$-hypergeometric term, we associate it with a
normal representation ($q$-NR) which consists of four polynomials
$r,s,u,v$. Based on the properties of the representation, we may
give a definition of $q$-proper hypergeometric terms and prove
that under the condition that $v$ is $\e_y$-free, a bivariate
$q$-hypergeometric term has a $qZ$-pair if and only if it is a
$q$-proper term. Applying the decomposition \eqref{ddd}, we deduce
that for any bivariate $q$-hypergeometric term $T$, it has a
$qZ$-pair if and only if $T_2$ is $q$-proper.

We conclude with some examples.

\section{$\e$-Free Decomposition}
\label{q-decom} Throughout the paper, we let $\Z,\Z^+$ and $\N$
denote the set of integers, positive integers and nonnegative
integers, respectively. For integers (or polynomials) $a,b$, we
denote by $\gcd(a,b)$ the (monic) greatest common divisor of $a$
and $b$. We also write $a \, \bot \, b$ to indicate that $a$ and
$b$ are relatively prime, i.e., $\gcd(a,b)=1$.

Let $\F$ be a field of characteristic zero, $q \in \F$ a nonzero
element which is not a root of unity, and $x$ transcendental over
$\F$. Denote by $\e$ the unique automorphism of $\F(x)$ which
fixes $\F$ and satisfies $\e x=qx$. Then $\F(x)$ together with the
{\it $q$-shift operator} $\e$ is a difference field \cite{Cohn}.
Let $r$ and $s$ be two polynomials. We say that $r/s$ is {\it
$\e$-reduced} if $r \,\bot\, \e^h s$ for all $h \in \Z$.

To be more specific, the rational functions involved in the
$q$-hypergeometric terms (see Definition \ref{hypdef}) are
rational functions of $q^n$. However, for a rational function $R
\in \F(x)$, we have
$$N\,R(q^n)=R(q^{n+1}) = \e R(q^n)\quad \mbox{and} \quad  R(q^n)=0 \ \forall\, n \ge n_0
\Leftrightarrow R(x)=0.$$ Therefore, there is a natural one-to-one
correspondence between the set of rational functions of $q^n$
together with the shift operator $N$ and the field $\F(x)$
together with the $q$-shift operator $\e$. In this paper, we adopt
the notation of $\F(x)$ as in the work of
Abramov-Paule-Petkov\v{s}ek \cite{Abra-P-P}.

The concept of rational normal forms introduced by Abramov and Petkov\v{s}ek
\cite{Ab-Pet} can be extended to the $q$-case.
\begin{defi}
\label{def-qRNF} Let $R \in \F(x)$ be a rational function. If polynomials
$r,s,u,v \in \F[x]$ satisfy
\begin{itemize}
\item[{\rm (i)}] $R={r \over s} \cdot {\e (u/v) \over
(u/v)}$ where $u \,\bot\, v$ and $u,v$ have no factor $x$,
\item[{\rm (ii)}] $r/s$ is $\e$-reduced,
\end{itemize}
then $(r,s,u,v)$ is called a $q$-rational normal form {\rm
(}$q$-RNF{\rm )} of $R$.
\end{defi}

Recall that a monic polynomial that has no factor $x$ is called a
$q$-monic polynomial by  Abramov, Paule, and Petkov\v{s}ek
\cite{Abra-P-P}. The following factorization theorem was given in
\cite{Abra-P-P}.

\begin{theo}
\label{fac-R} Let $R \in \F(x) \setminus \{0\}$. Then there exist
$z \in \F$ and monic polynomials $a,b,c \in \F[x]$ such that
\begin{equation}
\label{3-condi}
\begin{array}{l}
 R(x) = \displaystyle{z {a(x) \over b(x)} {c(qx) \over c(x)}},\\[8pt]
 \gcd(a(x),b(q^nx))=1, \quad \mbox{for all $n \in \N$},\\[5pt]
 \gcd(a(x),c(x))=\gcd(b(x),c(qx))=1 \quad \mbox{and} \quad  c(0) \not = 0.
\end{array}
\end{equation}
\end{theo}
We call $(az,b,c)$ a {\it $q$-Gosper form} ($q$-GF) of $R$.

\begin{theo}
\label{exist-RNF}
Every rational function $R \in \F(x)$ has a $q$-RNF.
\end{theo}
\pf It is clear that $(0,1,1,1)$ is a $q$-RNF of $0$. For
$R\not=0$, by Theorem~\ref{fac-R}, there exists a $q$-GF
$(az,b,c)$ of $R$. Applying Theorem~\ref{fac-R} again to
$b(x)/a(x)$, we get a $q$-GF $(r,s,d)$. From the construction
given in \cite{Abra-P-P}, we have $r\,|\,b$ and $s\,|\,a$. Hence
$s(x) \,\bot\, r(xq^n)$ for any $n \in \N$ because $(az,b,c)$ is a
$q$-GF. Since $(r,s,d)$ is also a $q$-GF, we have $r(x) \,\bot\,
s(xq^n)$ for any $n \in \N$. Thus $s/r$ is $\e$-reduced and
$(zs,r,c/\gcd(c,d),d/\gcd(c,d))$ is a $q$-RNF of $R$. \qed

The above proof provides an algorithm to generate a $q$-RNF of
$R$.

{\noindent \bf Algorithm $q$-RNF} \\
\indent if $R=0$ then \\
\indent \indent  return $(0,1,1,1)$; \\
\indent else \\
\indent \indent compute `$q$-GF' of $R$, we get $(a,b,c)$; \\
\indent \indent compute `$q$-GF' of $b/a$, we get $(r,s,d)$; \\
\indent \indent return $(s,r,c/\gcd(c,d),d/\gcd(c,d))$; \\

We now come to the  $q$-multiplicative representation of a general
$q$-hypergeometric term. This is the starting point of the
$\e$-free decomposition algorithm.

\begin{defi}
\label{hypdef} Suppose $T(n)$ is a function from $\N$ to $\F$. If there exist a
nonnegative integer $n_0$ and a nonzero rational function $R(x) \in \F(x)$ such
that $T(n+1)=R(q^n)T(n)$ for all $n \ge n_0$, then we call $T(n)$ a
{\rm(}univariate\rm{)} $q$-hypergeometric term.
\end{defi}

Suppose $(r,s,u,v)$ is a $q$-RNF of a rational function $R$.
Then the corresponding $q$-hypergeometric
term $T(n)$ satisfies
$$
T(n)=T(n_0)\prod_{j=n_0}^{n-1} R(q^j)= {T(n_0) \over u(q^{n_0})/v(q^{n_0})} \cdot
   {u(q^n) \over v(q^n)} \prod_{j=n_0}^{n-1} {r(q^j) \over s(q^j)}, \ \forall\, n \ge n_0.
$$
This leads to the following definition.
\begin{defi}
Let $T(n)$ be a $q$-hypergeometric term and $D,U$ be two rational
functions such that $D(q^n)$ has neither poles nor zeros and
$U(q^n)$ has no poles for all $n \ge n_0$. Suppose that
 \[
T(n)=U(q^n) \prod_{j=n_0}^{n-1} D(q^j), \qquad \forall\, n \ge
n_0.\]
Then we call $(D,U,n_0)$ a $q$-multiplicative
representation {\rm(}$q$-MR{\rm)} of $T$.
\end{defi}

Let $\Delta=N-1$ be the difference operator with respect to $n$.
The following lemma can be easily verified.

\begin{lem}
\label{q-trans} Let $T$ and $T_1$ be two $q$-hypergeometric terms
with $q$-MRs $(D,U,n_0)$ and $(D,U_1,n_0)$, respectively. Suppose
that \[ T_2=T- \Delta T_1\quad \mbox{ and } \quad U_2=U-D \cdot \e
U_1 + U_1. \]
 Then $(D,U_2,n_0)$ is a $q$-MR of $T_2$.
\end{lem}

For $u,v \in \F[x]$, let ${\cal R}$ be the set
of all nonnegative integers $h$ such that there exists an irreducible
polynomial $p(x) \not = x$ satisfying
$p(x)\mid u(x)$ and $p(x)\mid v(q^h x)$.
Define $\qdis(u,v)$ to be $\max \{h \in {\cal R} \}$ or $-1$ if $\cal R$
is empty. Note that ${\cal R}$ is a finite set, and ``$\qdis$'' is well defined.
If $\qdis(v,v)=0$, we say that $v$ is {\it $\e$-free}.

Given a $q$-hypergeometric term $T$ with a $q$-MR $(D,U,n_0)$.
Usually the denominator $u$ of $U$ is not $\e$-free. However,
translating the decomposition algorithm of \cite{Ab-Pet} into the
$q$-case, we have the following $\e$-free decomposition algorithm
``$q$-decomp'', which decomposes $T=\Delta T_1+T_2$ such that
$T_2$ has a $q$-MR $(F,V,n_0)$ where the denominator of $V$ is
$\e$-free.

\noindent {\bf Algorithm $q$-decomp}

\noindent Input: $(D,U,n_0)$  \qquad  Output: $U_1,F,V \in \F(x)$

\vskip 5pt

\noindent
\indent $d_1:=\num(D)$; $d_2:=\den(D)$; \\
\indent $U_1:=0$; $U_2:=U$; $u_2:=\den (U)$; \\
\indent $N:=\qdis(u_2,u_2)$;\\
\indent for $h:=N$ down to $1$ do \\
\indent \indent $v_2:=u_2/\gcd(u_2,d_2)$; \\
\indent \indent $s(x):=\gcd(v_2(x),v_2(q^{-h}x))$; \\
\indent \indent $(\tilde{s},\tilde{u}_2) := \mbox{pump}(s,u_2)$; \\
\indent \indent write $U_2=a/\tilde{u}_2 + b/\tilde{s}$
            where $a, b \in \F[x]$; \\
\indent \indent $U_1':=-b/\tilde{s}$;\\
\indent \indent $U_1:=U_1+U_1'$; $U_2:=U_2-D \cdot \e U_1' + U_1'$; \\
\indent \indent $u_2:=\den(U_2)$; \\
\indent $f_1:=d_1$; $f_2:=d_2$; $v_1:=\num(U_2)$; $v_2:=\den(U_2)$; \\
\indent $w:=\gcd(d_2,v_2)$;\\
\indent $v_2:=v_2/w$; $f_2:=\e w f_2/w$;\\
\indent $F:=f_1/f_2$; $V:=(1/w(q^{n_0})) \cdot v_1/v_2$; \\
\indent return $(U_1,F,V)$.

\noindent The procedure ``pump'' is the same as in the ordinary case.

\noindent {\bf Algorithm pump} \\
Input: $f,g \in \F[x]$; \qquad Output: $\tilde{f},\tilde{g} \in \F[x]$.

\vskip 5pt

\indent $\tilde{f}:=f; \tilde{g}:=g/f;$ \\
\indent repeat \\
\indent \indent $d:=\gcd(\tilde{f}, \tilde{g});$
   \quad $\tilde{f}:=\tilde{f}d; \tilde{g}:=\tilde{g}/d;$ \\
\indent until $\deg d= 0$; \\
\indent return ($\tilde{f},\tilde{g}$).

The following theorem shows that the $\e$-free algorithm generates
the desired decomposition.

\begin{theo}
\label{q-decomp}
Let $T$ be a $q$-hypergeometric term with a $q$-MR $(D,U,n_0)$
and $U_1,F,V$ be given by the algorithm $q$-decomp.
Then there exist
$q$-hypergeometric terms $T_1$ and $T_2$ such that
\begin{itemize}
\item[{\rm (1)}]
$T= \Delta T_1 + T_2$.
\item[{\rm (2)}]
$T_1$ has a $q$-MR $(D,U_1,n_0)$
and $T_2$ has a $q$-MR $(F,V,n_0)$.
\item[{\rm (3)}]
The denominator of $V$ is $\e$-free.
\end{itemize}
Furthermore, if $D$ is $\e$-reduced, so is $F$.
\end{theo}
\pf
Let $u_0$ be the denominator of $U$.
We first use induction to show that after iterating the loop of $h$ in the algorithm
$i$ times,
the denominator $u_2$ of $U_2$ satisfies:
\begin{itemize}
\item[(a)] $\qdis(v_2,v_2) \le N-i$,
\item[(b)] $u_2(q^n)$ has no zeros for all $n \ge n_0$,
\end{itemize}
where $v_2=u_2/\gcd(u_2,d_2)$, and $d_2$ is the denominator of $D$.

The case for $i=0$ is trivial. Assume that the assertion holds for
$i-1$. Let $u_2$ and $u'_2$ be the denominator of $U_2$ after
$i-1$ and $i$ iterations, respectively. Set $h=N-(i-1)>0$ and
$w_2=\gcd(u_2,d_2)$. From the algorithm $q$-decomp we have
\[ v_2=u_2/w_2 \quad \mbox{ and } \quad s=\gcd(v_2(x), v_2(q^{-h}x)).\]
 Suppose the
prime decomposition of $s$ is $p_1^{\alpha_1} \cdots
p_r^{\alpha_r}$ and $v_2=p_1^{\beta_1} \cdots p_r^{\beta_r} v',
w_2=p_1^{\gamma_1} \cdots p_r^{\gamma_r} w'$ where $v' \,\bot\, s,
w' \,\bot\, s$. Then the algorithm ``pump'' enables  us to
decompose $u_2$ as $p_1^{\beta_1+\gamma_1} \cdots
p_r^{\beta_r+\gamma_r} \cdot (v'w')$. That is,
$\tilde{s}=p_1^{\beta_1+\gamma_1} \cdots p_r^{\beta_r+\gamma_r}$
and $\tilde{u}_2=v'w'$. Since \[ U_2={a \over \tilde{u}_2} + {d_1
\over d_2} \cdot \e \left( {b \over \tilde{s}} \right),\]
 it follows
that $u'_2$ divides the least common multiple of $\tilde{u}_2$ and
$d_2 \e \tilde{s}$. Hence we have that  $u'_2$ divides $v' d_2
\cdot \e \tilde{s}$. Let $v''=v' \cdot \e \tilde{s}$. Assume that
there exist an integer $m \ge h$ and an irreducible polynomial
$p(x) \not= x$ such that $p \mid v''$ and $p \mid \e^m v''$. We
may encounter four cases:
\begin{itemize}
\item
$p \mid v'$ and $p \mid \e^m v'$. \\
From $v' \mid v_2$ and $\qdis(v_2,v_2) \le h$, it follows that $m=h$.
Therefore, $\e^{-h} p \mid \e^{-h} v_2$ and $\e^{-h} p \mid v_2$. Consequently,
we have $\e^{-h} p \mid s$, which contradicts $v' \,\bot\, s$. \item
$p \mid v'$ and $p \mid \e^{m+1} \tilde{s}$. \\
Since $s$ and $\tilde{s}$ have the same prime factors, we have $p \mid \e^{m+1}
s$, implying that  $p \mid \e^{m+1} v_2$. On the other hand, we have $p \mid
v_2$, which contradicts $\qdis(v_2,v_2) \le h$. \item
$p \mid \e \tilde{s}$ and $p \mid \e^m v'$. \\
In this situation, we have $\e^{-1} p \mid \tilde{s}$, which implies that
$\e^{-1} p \mid \e^{-h} v_2$, or equivalently, $\e^{h-1} p \mid v_2$. On the
other hand, $\e^{h-1} p \mid \e^{m+h-1} v_2$. Since $\qdis(v_2,v_2) \le h$, we
get $m+h-1 \le h$,  and hence $m=1$. Now we have $p \mid \e s$ and $p \mid \e
v'$, which contradicts $v' \,\bot\, s$. \item
$p \mid \e \tilde{s}$ and $p \mid \e^{m+1} \tilde{s}$. \\
Similarly, we have $\e^{-1} p \mid s$ and hence $\e^{-1} p \mid \e^{-h} v_2$,
i.e., $\e^{h-1} p \mid v_2$. However, we have $\e^{h-1} p \mid \e^{m+h} v_2$.
Thus, we obtain $m+h \le h$,  which is also a contradiction.
\end{itemize}
In summary, we may conclude that $\qdis(v'',v'') \le h-1$. Because
$u'_2$ divides $v'' \cdot d_2$,  there exist $\bar{v} \mid v''$
and $\bar{w} \mid d_2$ such that $u'_2=\bar{v} \bar{w}$. Let
$v_2'=u'_2/\gcd(u'_2,d_2)$. From $\bar{w} \mid \gcd(u'_2,d_2)$, it
follows that $v'_2 \mid  \bar{v}$. So we get $\qdis(v_2',v_2') \le
h-1=N-i$. Thus, we have proved (a). Since $u_2' | u_2 \cdot \e u_2
\cdot d_2$, (b) immediately follows from the induction hypothesis.

On the other hand, since  $\tilde{s} \mid  u_2$, (b) implies that $U_1(q^n)$ has no
poles for all $n \ge n_0$. Let
\begin{equation}
\label{t1-t2} T_1(n) = U_1(q^n) \prod_{j=n_0}^{n-1} D(q^j) \quad \mbox{and}
\quad T_2(n) = U_2(q^n) \prod_{j=n_0}^{n-1} D(q^j).
\end{equation}
Noting that $U_2=U-D \e U_1 + U_1$, by Lemma~\ref{q-trans}, we obtain $T=\Delta
T_1+T_2$.

Because $w \mid  d_2$ and $d_2(q^n) \not = 0$ for all $n \ge n_0$, we can
write $T_2(n)$ as
$$
T_2(n) = {1 \over w(q^{n_0})} U_2(q^n) w(q^n)
 \prod_{j=n_0}^{n-1} D(q^j) {w(q^j) \over w(q^{j+1})}
  = V(q^n) \prod_{j=n_0}^{n-1} F(q^j).
$$

Let $v$ be the denominator of $V$. Then (a) implies $\qdis(v,v)=0$,
that is, $v$ is $\e$-free.

Finally, notice that $f_1=d_1$ and $f_2=\e w \cdot (d_2 / w)$ where $w \mid
d_2$. Therefore, $F$ is $\e$-reduced provided that $D$ is $\e$-reduced. This
completes the proof. \qed

\section{Bivariate $q$-Hypergeometric Terms}

We begin this section with the definition of bivariate $q$-hypergeometric
terms.

\begin{defi}
Suppose $T(n,k)$ is a function from $\N^2$ to
$\F$. If there exist
rational functions $R_1(x,y), R_2(x,y) \in \F(x,y)$ and
$n_0 \in \N$ such that
$$
T(n+1,k) = R_1(q^n,q^k) T(n,k) \quad \mbox{and} \quad
T(n,k+1) = R_2(q^n,q^k) T(n,k),
$$
for all $n,k \ge n_0$, then we call $T(n,k)$ a bivariate $q$-hypergeometric
term.
\end{defi}

Without loss of generality, from now on we may assume that $n_0=0$
and that $R_1(q^n,q^k), R_2(q^n,q^k)$ have neither zeros nor poles
for all $n,k \ge 0$.

Denote by $\e_x$ and $\e_y$ the shift operators on $\F(x,y)$ defined by $\e_x x
= qx$, $\e_x|_{\F(y)} = {\rm id}$ (the identity map) and $\e_y y=qy,
\e_y|_{\F(x)}={\rm id}$, respectively. The idea of $q$-RNF can be easily
adopted to the bivariate case by taking $\F(y)$ as the ground field. Let
$R(x,y)$ be a rational function of $x$ and $y$, its $q$-rational normal form
($q$-RNF with respect to $\e_x$) is represented by $(r,s,u,v)$ as in the
univariate case. By using the ground field $\F(x)$, we may find a $q$-RNF of
$R(x,y)$ with respect to $\e_y$.

Let $T(n,k)$ be a bivariate $q$-hypergeometric term. By definition, there
exists a rational function $R$ such that
$$T(n+1,k)/T(n,k)=R(q^n,q^k).$$
 Suppose $(r,s,u,v)$ is a $q$-RNF of $R$ with respect to $\e_x$. We
call $(r,s,u,v)$ a $q$-normal representation ($q$-NR) of $T(n,k)$ with respect
to the shift operator $N$. Similarly, we can define the $q$-NR of $T(n,k)$ with
respect to the shift operator $K$.

We next give a characterization of the polynomials involved in the
 $q$-NR of bivariate $q$-hypergeometric terms.

\begin{theo}
\label{main} Let $T(n,k)$ be a bivariate $q$-hypergeometric term that has a
$q$-NR $(r,s,u,v)$ with respect to $N$. Then $r$ and $s$ are products of
polynomials having the form
$$(x^c y^d) \cdot \prod_{l=1}^{a}  p (q^{w_l}x^{a}y^{b}),$$
where
$p$ is a Laurent polynomial of one variable, $a \in \Z^+, b,c,d, \, w_l \in \Z$,
$a \,\bot\, b$, and
$ w_i \not\equiv w_j \pmod{a},\ \forall \ i \not= j $.

Similarly, suppose $(r,s,u,v)$ is a $q$-NR of $T$ with respect to
$K$. Then $r$ and $s$ are products of polynomials having the form
\[  (x^c y^d) \cdot \prod_{l=1}^{a}  p (q^{w_l}x^{b}y^{a})\]
under the same conditions.
\end{theo}

{\noindent \it Sketch of the proof.}  The proof of the ordinary
case \cite[Theorem 3.4]{Hou} can be  carried over to the $q$-case
except that we need to consider the characterization of
polynomials $f(x,y)$ such that $f(q^ax,q^by) = C f(x,y)$ for
certain integers $a,b$ and $C \in \F$. \qed

Consequently, we have
\begin{coro}
\label{q-proper} Let $T(n,k)$ be a bivariate $q$-hypergeometric
term that has a $q$-NR $(r,s,u,v)$ with respect to $N$ {\rm(}or
$K$ respectively{\rm)}. Then we have
$$
T(n,k) = C \cdot {u(q^n,q^k) \over v(q^n,q^k)} \cdot {\prod\limits_{l=1}^{uu}
\prod\limits_{j=0}^{a_l n+b_l k+c_l} f_l(q^j)  \over
\prod\limits_{l=1}^{vv} \prod\limits_{j=0}^{a'_l n+b'_l k+c'_l} g_l(q^j)},
$$
where $C \in \F$, $uu, vv \in \N$,
$a_l,b_l,c_l,a'_l,b'_l,c'_l \in \Z$ and $f_l,g_l$ are polynomials.
\end{coro}

Corollary~\ref{q-proper} enables us to give the following
definition of $q$-proper hypergeometric terms.
\begin{defi}
A polynomial $f \in \F[x,y]$ is said to be $q$-proper if for each of its
irreducible factor $p(x,y) \in \F[x,y]$, there exist $a,b \in \Z$, not both
zeros, such that $p(x,y)|p(q^a x, q^b y)$. A bivariate $q$-hypergeometric term
$T$ is said to be $q$-proper if $v$ is a $q$-proper polynomial where
$(r,s,u,v)$ is a $q$-NR of $T$ with respect to $N$ or $K$.
\end{defi}
Suppose that $T$ is a bivariate $q$-hypergeometric term that has a
$q$-NR $(r,s,u,v)$ with respect to $N$ (or $K$).
Theorem~\ref{main} guarantees  that $r$ and $s$ are both
$q$-proper polynomials.

As in the case of ordinary bivariate hypergeometric terms (\cite[Theorems
4.2]{Hou}), we have  an analogous ``fundamental theorem'' for the $q$-case.
\begin{theo}
\label{fund} Let $T(n,k)$ be a bivariate $q$-hypergeometric term. Then $T$ is
$q$-proper if and only if there exist polynomials $a_{ij}(x) \in \F[x]$, not
all zero, such that
$$
\sum_{0 \le i \le I,\; 0 \le j \le J} a_{ij}(q^n) T(n+i,k+j) = 0
\quad \forall\, n,k \ge 0.
$$
\end{theo}

Based on an analogous argument for the ordinary case as in \cite[Theorem 6.2.1]{Pet-W-Z},
we get
\begin{coro}
\label{suffi} Any $q$-proper hypergeometric term has a $qZ$-pair.
\end{coro}

\section{The Existence of $qZ$-Pairs}

In this section, we obtain a necessary and sufficient condition for the
existence of $qZ$-pairs for any bivariate $q$-hypergeometric term based on its
$q$-NR with respect to $K$.

From Theorem~\ref{main}, we have

\begin{coro}
\label{coro-a/b} Let $T(n,k)$ be a bivariate $q$-hypergeometric term that has a
$q$-NR $(r,s,u,v)$ with respect to $K$. Then there exist polynomials $f_i(x), g_i(x)
\in \F[x]$ and $a_i,a'_i,b_i,b'_i \in \Z$ such that
$$
\prod_{j=0}^{k-1} \left( {r(q^{n+1},q^j) \over r(q^n,q^j)} \cdot {s(q^n,q^j)
\over s(q^{n+1}, q^j)} \right) = \prod_{i=1}^{\ell} {f_i(q^{a_i k + b_i n})
\over g_i(q^{a'_i k + b'_i n})}.
$$
\end{coro}

We need to consider the following ratio
$${T(n+i,k) \over T(n,k)}  =  {T(n+i,0) \over T(n,0)}
             \prod_{j=0}^{k-1} \left\{ {T(n+i,j+1)\over T(n+i,j)}
             {T(n,j)\over T(n,j+1)} \right\},
$$
which can be rewritten as
\begin{align}
\label{NT/T}
{T(n+i,k) \over T(n,k)} & = \prod_{l=0}^{i-1} \prod_{j=0}^{k-1}
\Big\{ {r(q^{n+l+1},q^j) \over  r(q^{n+l},q^j) }
    {s(q^{n+l},q^j) \over s(q^{n+l+1},q^j)} \Big\}
    \prod_{l=0}^{i-1}  {T(n+l+1,0) \over T(n+l,0)} \nonumber \\
& \qquad \cdot {u(q^{n+i},q^k) \over u(q^{n+i},q^0) } {u(q^n,q^0) \over u(q^n,q^k)}
{v(q^{n+i},q^0)\over v(q^{n+i},q^k)}  {v(q^n,q^k)\over v(q^n,q^0)}.
\end{align}

From  Corollary~\ref{coro-a/b} we get the following expression.

\begin{lem}
Let $T(n,k)$ be a bivariate $q$-hypergeometric term that has a $q$-NR
$(r,s,u,v)$ with respect to $K$. Then for each $i \ge 0$, there exist
$q$-proper polynomials $w_1^{(i)}(x,y)$ and $w_2^{(i)}(x,y)$ such that
\begin{equation}
\label{stand-def}
{T(n+i,k) \over T(n,k)} = {u(q^{n+i},q^k) \over v(q^{n+i},q^k)} \cdot
{v(q^n,q^k) \over u(q^n,q^k)} \cdot {w_1^{(i)}(q^n,q^k) \over w_2^{(i)}(q^n,q^k)},
\quad \forall\, n,k \ge 0.
\end{equation}
\end{lem}

An $\e_y$-free polynomial that is not $q$-proper has a special factor.

\begin{lem}
\label{d-Factor} Let $f \in \F[x,y]$ be a non-$q$-proper and
$\e_y$-free polynomial. Then there exists an irreducible
factor $p$ of $f$ such that
\begin{equation}
\label{decom}
\begin{array}{l}
p(x,y) \,\bot\, p(q^i x, q^j y), \quad \forall \, (i,j) \in \Z^2 \setminus \{(0,0)\}, \\
p(x,y) \,\bot\, f(q^i x, q^j y), \quad \forall \, (i,j)\in  (\N \times \Z) \setminus \{(0,0)\}.
\end{array}
\end{equation}
\end{lem}
\pf Since $f(x,y)$ is non-$q$-proper, by definition it has
an irreducible factor $p_1(x,y)$ such that $p_1(x,y) \, \bot \,
p_1(q^i x, q^j y),  \forall \, (i,j) \in \Z^2 \setminus \{(0,0)\}$.

We may factor $f(x,y)$ as
$$
f(x,y)=p_1^{\a_1}(q^{a_1}x,q^{b_1}y) \cdots p_1^{\a_r}(q^{a_r}x,q^{b_r}y)
f_1(x,y),
$$
where $(a_i,b_i) \in \Z^2$ are distinct pairs, $\alpha_i \in
\Z^+$, and $p_1(q^i x, q^j y) \, \bot \, f_1(x,y)$ for all $i,j
\in \Z$. Since $f(x,y)$ is $\e_y$-free, it follows that $a_i \neq
a_j$ as long as $i \neq j$. Without loss of generality, we may
assume that $a_1 < a_2 < \cdots <a_r$. Thus,
$p(x,y)=p_1(q^{a_1}x,q^{b_1}y)$ satisfies the condition
\eqref{decom}. \qed

We are now ready to give a criterion for the existence of
$qZ$-pairs.
\begin{theo}
\label{one-term} Let $T(n,k)$ be a bivariate $q$-hypergeometric term that has a
$q$-NR $(r,s,u,v)$ with respect to $K$ such that $v$ is $\e_y$-free. Then
$T(n,k)$ has a $qZ$-pair if and only if $v$ is a $q$-proper polynomial.
\end{theo}
\pf Because of Corollary~\ref{suffi}, it suffices to show that if
$T(n,k)$ has a $qZ$-pair, then it is $q$-proper. To this end, we
assume that $T(n,k)$ is a bivariate $q$-hypergeometric term.
Moreover, we assume that $T(n,k)$ is not $q$-proper, but it has a
$qZ$-pair. We proceed to find a contradiction.

Clearly, for a difference operator $L \in \F[q^n,N]$, we have
$$(N \cdot L) T(n,k)= (K-1) G(n,k) \Longleftrightarrow L T(n,k)=(K-1) G(n-1,k).$$
Therefore, we may assume that $T(n,k)$ has a $qZ$-pair $(L,G)$ of the form
\[L=\sum_{i=0}^I a_i(q^n) N^i,\]
 where $a_i(q^n)$ are polynomials in $q^n$
and $a_0 \not= 0$. Since $LT/T$ and $(K-1)G/G$ are both rational functions of
$q^n$ and $q^k$, we may assume that
$$G(n,k)= {f(q^n,q^k) \over g(q^n,q^k)} T(n,k) ,$$
where $f,g \in \F[x,y]$ are two relatively prime polynomials.

By the definition of $qZ$-pairs, we have
\begin{equation}
\label{$Z$-pair} \sum_{i=0}^I a_i(q^n) {T(n+i,k) \over T(n,k)} =
{f(q^n,q^{k+1}) \over g(q^n,q^{k+1})} {T(n,k+1) \over T(n,k)} - {f(q^n,q^k) \over
g(q^n,q^k)}.
\end{equation}
Substituting \eqref{stand-def} into
\eqref{$Z$-pair}, we obtain
\begin{equation}
\label{eq-1}
\sum_{i=0}^I a_i(x) {u(q^i x ,y) \over v(q^ix,y)}
{w_1^{(i)} (x,y) \over w_2^{(i)}(x,y)}
 = {f(x,qy) \over g(x,qy)} {r(x,y) \over s(x,y)} {u(x,qy) \over v(x,qy)}
 - {f(x,y) \over g(x,y)} {u(x,y) \over v(x,y)}.
\end{equation}

Let $u_1 = u / \gcd( u,g ),\,  g_1= g / \gcd(u,g)$.
Multiplying
$$
 g_1(x,qy) g_1(x,y) v(x,qy) s(x,y)
   \prod_{j=0}^I v(q^jx,y) w_2^{(j)}(x,y)
$$
to both sides of \eqref{eq-1}, we arrive at

\begin{align} \label{eq-3}
& g_1(x,qy) g_1(x,y) v(x,qy) s(x,y) \nonumber \\
& \quad \cdot \sum_{i=0}^I
     a_i(x) u(q^ix,y) w_1^{(i)}(x,y) \prod_{j \not= i}
           v(q^jx,y) w_2^{(j)}(x,y) \nonumber \\
= & f(x,qy) r(x,y) u_1(x,qy) g_1(x,y)
  \prod_{j=0}^I v(q^jx,y) w_2^{(j)}(x,y) \nonumber \\
& - f(x,y) u_1(x,y) g_1(x,qy) v(x,qy) s(x,y) w_2^{(0)}(x,y)
\cdot \prod_{j=1}^I v(q^jx,y)
w_2^{(j)}(x,y).
\end{align}

Since $T(n,k)$ is not $q$-proper,  from Lemma \ref{d-Factor} it
follows that there exists an irreducible factor $p$ of $v$
satisfying the condition \eqref{decom}. Noting that  $p(x,y)$
divides each term of the left-hand side of \eqref{eq-3} except for
the first term, we obtain that $p(x,y)$ divides
\begin{align*}
& g_1(x,qy) v(x,qy) s(x,y)
   \prod_{j=1}^I v(q^jx,y) w_2^{(j)}(x,y) \\
& \quad  \times \big( g_1(x,y)  a_0(x) u(x,y) w_1^{(0)}(x,y)
 + f(x,y) u_1(x,y) w_2^{(0)}(x,y) \big).
\end{align*}
From \eqref{decom} it follows that
 \[
p(x,y) \, \bot \, v(x,qy) \prod_{j=1}^I v(q^jx,y).\]
 Since $s$ and
$w_2^{(j)}$ are $q$-proper, they are also relatively prime to $p$.
This implies that $p(x,y)$ divides
\begin{align}
\label{r-A} g_1(x,qy) \big( g_1(x,y)  a_0(x) u(x,y) w_1^{(0)}(x,y)
 + f(x,y) u_1(x,y) w_2^{(0)}(x,y) \big).
\end{align}
Similarly, since $p(x,qy)$ divides both sides of
\eqref{eq-3} and $u \, \bot \, v$, we have
\begin{align}
\label{r-B} p(x,qy) \mid f(x,qy) g_1(x,y).
\end{align}
\begin{itemize}
\item[Case 1.] Suppose $p(x,qy) \mid f(x,qy)$. Since $p(x,y)$ divides
\eqref{r-A}, it follows that
$$
p(x,y) \mid  g_1(x,qy) g_1(x,y)  a_0(x) u(x,y) w_1^{(0)}(x,y).
$$
Since $f \, \bot \, g, u \, \bot \, v$,
$a_0$ and $w_1^{(0)}$ are $q$-proper polynomials, we may deduce
that $p(x,y) \mid g_1(x,qy)$, i.e., $p(x,q^{-1}y) \mid g_1(x,y)$. Let
$m(>\!0)$ be the greatest integer such that $p(x,q^{-m}y) \mid g_1(x,y)$.
By virtue of \eqref{eq-3}, we have that $p(x,q^{-m}y)$ divides
\begin{align*}
f(x,y) u_1(x,y) g_1(x,qy) v(x,qy) s(x,y) w_2^{(0)}(x,y)
  \prod_{j=1}^I v(q^jx,y) w_2^{(j)}(x,y).
\end{align*}
However, $f \, \bot \, g$ and $g_1 \, \bot \,
u_1$ imply that $p(x,q^{-m}y) \mid g_1(x,qy)$, which contradicts
the choice of $m$.

\item[Case 2.] Suppose $p(x,qy) \mid g_1(x,y)$. Let $M>0$ be the
greatest integer such that $p(x,q^My) \mid g_1(x,y)$. Similarly,
from \eqref{eq-3} it follows that
$p(x,q^{M+1}y)$ divides
$$
f(x,qy) r(x,y) u_1(x,qy) g_1(x,y)
  \prod_{j=0}^I v(q^jx,y) w_2^{(j)}(x,y).
$$
Hence we get $p(x,q^{M+1}y) \mid g_1(x,y)$, which is again a
contradiction. \qed
\end{itemize}

To extend the above result to general bivariate $q$-hypergeometric
terms, we need the concept of similar $q$-hypergeometric terms.
Two bivariate $q$-hypergeometric terms $T_1,T_2$ are called {\it
similar} if there exists a rational function $R \in \F(x,y)$ such
that $T_1(n,k)/T_2(n,k)=R(q^n,q^k)$.

As in the ordinary case, the existence of $qZ$-pairs is preserved
under addition of similar bivariate $q$-hypergeometric terms.
\begin{lem}
\label{add} Suppose there exist $qZ$-pairs for two similar
bivariate $q$-hypergeometric terms $T_1(n,k)$ and $T_2(n,k)$. Then
there exists a $qZ$-pair for $T(n,k)=T_1(n,k) + T_2(n,k)$.
\end{lem}

Notice that $T(n,k)=(K-1)G(n,k)$ has a $qZ$-pair $(1,G)$. Combining
Theorem~\ref{one-term} and Lemma~\ref{add}, we obtain the main result of this
paper.
\begin{theo}
\label{criterion} Let $T(n,k)$ be a bivariate $q$-hypergeometric term. Let
$T_1,T_2$ be two similar bivariate $q$-hypergeometric terms satisfying
$$
T(n,k) = (K-1) T_1(n,k) + T_2(n,k)
$$
and $T_2(n,k)$ has a $q$-NR $(r,s,u,v)$ with respect to $K$ such that $v$ is
$\e_y$-free. Then $T(n,k)$ has a $qZ$-pair if and only if $T_2(n,k)$ is a
$q$-proper hypergeometric term, or equivalently, if and only if $v(x,y)$ is a
$q$-proper polynomial.
\end{theo}

\section{Algorithms}
\label{Alg} Let $T(n,k)$ be a bivariate $q$-hypergeometric term. By the
algorithm ``$q$-RNF'', we may find a $q$-NR $(r,s,u,v)$ of $T(n,k)$ with
respect to $K$. Let
$$
F(k) = {u(x,q^k)
\over v(x,q^k)}
  \prod_{j=0}^{k-1} {r(x,q^j) \over s(x,q^j)}, \quad \forall\, k \in \N.
$$
Then $F(k)$ is a univariate $q$-hypergeometric term over the field $\F(x)$ with
a $q$-MR $(r/s,u/v,0)$. On the other hand, by Equation~\eqref{NT/T}, we have
\begin{multline*}
{F(k)|_{x=q^{n+1}} \over F(k)|_{x=q^n}} = {u(q^{n+1},q^k) v(q^n,q^k) \over
u(q^n,q^k) v(q^{n+1},q^k)} \prod_{j=0}^{k-1} {r(q^{n+1},q^j) s(q^n,q^j) \over
r(q^n,q^j) s(q^{n+1},q^j)} \\
= {T(n+1,k) \over T(n,k)} \cdot {T(n,0) \over T(n+1,0)} \cdot {u(q^{n+1},q^0)
v(q^n,q^0) \over u(q^n,q^0) v(q^{n+1},q^0)},
\end{multline*}
which is also a rational function on $q^n$ and $q^k$. Hence $\widetilde{F}(n,k)
= F(k)|_{x=q^n}$ is a bivariate $q$-hypergeometric term.

Using the algorithm ``$q$-decomp'' given in Section~\ref{q-decom},
one may  find univariate $q$-hypergeometric terms $F_1(k), F_2(k)$
such that
$$
F(k) = (K-1) F_1(k) + F_2(k)
$$
and $F_2(k)$ has a $q$-MR $(f_1/f_2,v_1/v_2,0)$ with $v_2$ being
$\e_y$-free. Since $f_1/f_2, v_1/v_2 \in \F(x)(y)$, we may assume
that $f_1,f_2,v_1,v_2 \in \F[x,y]$ and $f_1 \,\bot\, f_2, v_1
\,\bot\, v_2$. From the fact that $r/s$ is $\e_y$-reduced, it
follows that $f_1/f_2$ is also $\e_y$-reduced.

Let
\begin{align*}
T_1(n,k)=T(n,0) {v(q^n,q^0) \over u(q^n,q^0)} \cdot F_1(k)|_{x=q^n}
, \\
T_2(n,k)=
T(n,0) {v(q^n,q^0) \over u(q^n,q^0)}
\cdot F_2(k)|_{x=q^n}.
\end{align*}
Since Equation~\eqref{t1-t2} implies that
\[
F_1(k) = {U_1 \over u/v} \cdot F(k) \quad \mbox{and} \quad F_2(k) = {v_1/v_2
\over u/v} \cdot F(k),
\]
it follows that $T_1(n,k)$ and $T_2(n,k)$ are similar bivariate
$q$-hypergeometric terms. It is easily verified that
$$
T(n,k) = (K-1) T_1(n,k) + T_2(n,k)
$$
and $(f_1,f_2,v_1,v_2)$ is a $q$-NR of $T_2$ with respect to $K$. Therefore,
Theorem~\ref{criterion} implies that $T(n,k)$ has a $qZ$-pair if and only if
$v_2$ is a $q$-proper polynomial.

Finally, we need the algorithm given by Le \cite{Le}
for determining whether or not a
polynomial is $q$-proper.

We are now ready to describe the algorithm to determine whether a bivariate
$q$-hypergeometric term $T(n,k)$ has a $qZ$-pair. {\tt
\begin{itemize}
\item[1.]
Apply the algorithm in \cite{Bo-Ko} to find a rational function $R
\in \F(x,y)$ such that
$${T(n,k+1) \over T(n,k)} = R(q^n,q^k).$$
\item[2.]
Find a $q$-RNF $(r,s,u,v)$ with respect to $\e_y$ of $R$.
\item[3.]
For $D=r/s, U=u/v$ and $n_0=0$, apply the algorithm `$q$-decomp' with respect
to $\e_y$ to get $V=v_1/v_2$.
\item[4.]
Use the algorithm in \cite{Le} to
determine whether $v_2$ is $q$-proper. If the answer is yes, then $T$ has a
$qZ$-pair; otherwise, $T$ does not have any $qZ$-pair.
\end{itemize}
}

Here are two examples.

\noindent {\bf Example 1.} Let
$$
T(n,k) = {q^k(1+q^{n+1}+q^{k+2}) \over (q^n+q^k+1)(q^n+q^{k+1}+1)\prod_{j=1}^{k+1} (1-q^j)}.
$$
Then
$$
{T(n,k+1) \over T(n,k)} =
{q (1+q^{n+1}+q^{k+3})(q^n+q^k+1) \over (q^n+q^{k+2}+1)(1+q^{n+1}+q^{k+2})(1-q^{k+2})},
$$
and we have
$$
r=q, \ s=1-q^2y, \ u=1+qx+q^2y, \ v=(x+y+1)(x+qy+1)
$$
is a $q$-NR of $T$ with respect to $K$. For $D=r/s, U=u/v$ and $n_0=0$,
applying the algorithm ``$q$-decomp'', we get
$$V=v_1/v_2 = {-q^2 \over (-1+q^2)(x+1)}.$$
Clearly,  $v_2$ is  $q$-proper, so
$T(n,k)$ has a $qZ$-pair. Indeed, we can check that
$$ L = 1, \quad G = {1 \over (q^n+q^k+1) \prod_{j=1}^k (1-q^j)}$$
is a $qZ$-pair for $T(n,k)$.

\noindent {\bf Example 2.} Let
$$
T(n,k) = {q^k(1+q^{n+1}+q^{k+2}) \over (q^n+q^k+1)(q^n+q^{k+1}+1)\prod_{j=1}^k (1-q^j)}.
$$
Then
$$
{T(n,k+1) \over T(n,k)} =
{q (1+q^{n+1}+q^{k+3})(q^n+q^k+1) \over (q^n+q^{k+2}+1)(1+q^{n+1}+q^{k+2})(1-q^{k+1})},
$$
and we have
$$
r=q, \ s=1-qy, \ u=1+qx+q^2y, \ v=(x+y+1)(x+qy+1)
$$
is a $q$-NR of $T$ with respect to $K$. For $D=r/s, U=u/v$ and $n_0=0$,
applying the algorithm ``$q$-decomp'', we get
$$V=v_1/v_2 = {-(x+y+1)q^2 \over (q-1)(x+1)(x+qy+1)}.$$
Since $x+qy+1$ is not a $q$-proper polynomial, it follows that
$T(n,k)$ has no $qZ$-pair.

\vskip 15pt {\small {\bf Acknowledgments.} The authors are
grateful to the referees for many helpful comments and
suggestions. We thank Professor Marko Petkov\v{s}ek for
bringing our attention to the work of Ore and Sato. This work was done under the
auspices of the National ``973'' Project on Mathematical
Mechanization, the National Science Foundation, the Ministry of
Education, and the Ministry of Science and Technology of China.}

\end{document}